\documentclass[secthm,seceqn,amsthm,ussrhead]{amsart}
\usepackage{amsmath,latexsym}
\usepackage[psamsfonts]{amssymb}
\usepackage{times}
\usepackage[mathcal]{euscript}
\numberwithin{equation}{section}

\newcommand{\bbT}{\mathbb T}

\renewcommand{\epsilon}{\varepsilon}

\newcommand{\be}{\begin{equation}}
\newcommand{\ee}{\end{equation}}
\newcommand{\no}{\nonumber}

\newcommand{\C}{\mathbb{C}}

\newcommand{\F}{\mathbb{F}}

\newcommand{\N}{\mathbb{N}}

\newcommand{\R}{\mathbb{R}}
\renewcommand{\S}{\mathbb{S}}
\newcommand{\T}{\mathbb{T}}

\newcommand{\Z}{\mathbb{Z}}

\newcommand{\cE}{{\mathcal E}}
\newcommand{\cF}{{\mathcal F}}

\newcommand{\cH}{{\mathcal H}}

\newcommand{\cU}{{\mathcal U}}

\renewcommand{\Re}{{\ensuremath{\mathrm{Re}}}}

{\bf}{\it}
\newtheorem{theorem}{Theorem}[section]
\newtheorem{lemma}[theorem]{Lemma}
\newtheorem{corollary}[theorem]{Corollary}

\newtheorem{definition}[theorem]{Definition}

\newtheorem{remark}[theorem]{Remark}

\date{\today}

\begin{document}
\title[Asymptotics for the number of eigenvalues of three-particle
Schr\"{o}dinger...]
 {Asymptotics for the number of eigenvalues of three-particle
  Schr\"{o}dinger operators on lattices}


 \author{Sergio Albeverio$^{1,2,3}$, G.F.Dell Antonio$^{4}$,
 Saidakhmat
 N. Lakaev$^{5,6}$}

\address{$^1$ Institut f\"{u}r Angewandte Mathematik,
Universit\"{a}t Bonn (Germany)} \email{albeverio@uni.bonn.de}
\address{
$^2$ \ SFB 611, \ Bonn, \ BiBoS, Bielefeld - Bonn\ (Germany)}
\address{
$^3$  CERFIM, Locarno and USI (Switzerland)}

\address{
{$^4$ Dipartment di Matematica, Univ.of Roma 1 and SISSA (Trieste)}}

\address{
$^5$ Samarkand Division of Academy of sciences of Uzbekistan
(Uzbekistan)} \email{lakaev@yahoo.com}

\address{
$^6$ {Samarkand State University,Samarkand (Uzbekistan)} \ {E-mail:
slakaev@mail.ru }}

\begin{abstract}
We consider the Hamiltonian of a system of three quantum mechanical
particles (two identical fermions and boson)on the three-dimensional
lattice $\Z^3$ and interacting by means of zero-range attractive
potentials. We describe the location and structure of the essential
spectrum of the three-particle discrete Schr\"{o}dinger operator
$H_{\gamma }(K),$ $K$ being the total quasi-momentum and $\gamma>0$
the ratio of the mass of fermion and boson.

 We choose  for $\gamma>0$ the interaction $v(\gamma)$ in such a
 way  the system consisting
 of one fermion and one boson  has a zero energy resonance.

We prove for any $\gamma> 0$ the existence infinitely many
  eigenvalues of the operator $H_{\gamma }(0).$

We establish for the number $N(0,\gamma; z; )$ of eigenvalues lying
below $z<0$ the following asymptotics
$$
\lim_{z\rightarrow 0-}\frac{N(0,\gamma;z )}{\mid \log \mid z\mid
\mid }={\ U} (\gamma ) .$$

Moreover, for all  nonzero values of the quasi-momentum $K \in T^3 $
we establish the finiteness of the number $
N(K,\gamma;\tau_{ess}(K))$ of eigenvalues of $H(K)$ below the bottom
of the essential spectrum
 and we give an asymptotics for the number $N(K,\gamma;0)$
of eigenvalues  below zero.
\end{abstract}

\maketitle



\section{Introduction}
We consider a system of three particles (two identical fermions and
a boson) on the three-dimensional lattice $\Z^3$ interacting by
means of zero-range attractive potentials.

The main goal of the present paper is to prove (under the conditions
relevant for Evimov's effect) the finiteness or infiniteness of the
number of eigenvalues lying below the bottom of the essential
spectrum of the three-particle discrete Schr\"{o}dinger operator
$H_{\gamma }(K)$ depending on the total quasi-momentum $K\in\T^3$
and the ratio $\gamma>0$ of the mass of the fermions and the boson.

 Efimov`s  effect is one of the remarkable results in the
spectral analysis for continuous and discrete  three-particle
Schr\"{o}dinger operators: if none of the three two-particle
Schr\"{o}dinger operators (corresponding to the two-particle
subsystems) has negative eigenvalues, but at least two of them have
a zero energy resonance, then this three-particle Schr\"{o}dinger
operator has an infinite number of discrete eigenvalues,
accumulating at zero.

Since its discovery by Efimov in 1970 \cite{Efi} much research have
been devoted to this subject. See, for example
\cite{AHW,AmNo,DFT,FaMe,OvSi,Sob,Tam91,Tam94,Yaf74}.

The main result obtained by  Sobolev \cite{Sob} (see also \cite{Tam
94}) is an asymptotics of the form $\cU_0|log|\lambda||$ for the
number of eigenvalues below $\lambda,\lambda<0$, where the
coefficient ${\cU}_0$ does not depend on the two-particle potentials
$ v_\alpha $ and is a positive function of the ratios
$m_1/m_2,m_2/m_3$ of the masses of the three-particles.

Recently  the existence of Efimov's effect for $N$-body quantum
systems with $N\geq 4$ has been proved by X.P. Wang in \cite{Wang}.

In fact in \cite{Wang} for the total (reduced) Hamiltonian a lower
bound on the number of eigenvalues of the form
$$C_0|log(E_0-\lambda)|$$ is given, when $\lambda$ tends to $E_0$,
where $C_0$ is a positive constant and $E_0$ is the bottom of the
essential spectrum.

 The kinematics of quantum particles
on lattices, even in the two and three-particle sector, is rather
exotic. For instance,
 due to the fact that   the discrete analogue of
the Laplacian or its generalizations (see \eqref{free}) are not
translationally invariant, the Hamiltonian of a system  does not
separate into two parts, one relating to the center-of-mass motion
and the other one relating to the internal degrees of freedom.

As a consequence any local substitute of the effective mass-tensor
(of a ground state) depends on the quasi-momentum of the system and,
in addition, it  is only semi-additive (with respect to the partial
order on the set of positive definite matrices). This is the
so-called {\it excess mass}
 phenomenon for lattice systems
(see, e.g., \cite{Mat} and  \cite{Mog91}):
 the effective mass of the bound state
of an $N$-particle system is greater than and in general,
 not equal to the sum of the effective
masses of the constituent quasi-particles.

The three-body  problem on lattices  can  be reduced to
three-particle Schr\"odinger operators  by using the Gelfand
transform. The underlying Hilbert space $\ell^2((\Z^d)^3)$
 is decomposed as a direct von Neumann integral associated with the
representation of the discrete group $\Z^3$  by shift operators on
the lattice and  the total three-body Hamiltonian appears to be
decomposable. In contrast to the continuous case, the corresponding
fiber Hamiltonians $H(K) $ associated with  the direct decomposition
depend parametrically  the quasi-momentum,  $K\in
\T^3=(-\pi,\pi]^3$, which  ranges over a cell of the dual lattice
(see \cite{RSIV}). Due to the loss of the spherical symmetry of the
problem, the spectra of the family $H(K)$ turn out to be  rather
sensitive to the quasi-momentum $K$.

In particular, Efimov's effect exists only for the zero value of the
three-particle quasi-momentum $K$ (see, e.g.,
\cite{ALM98,ALzM04,Ltmf91,Lfa93, LaAb99,Mat} for relevant
discussions and \cite{ALMM,FIC,GrSc,KM,Mat,MS,Mog91,RSIII} for the
general study of the low-lying excitation spectrum for quantum
systems on lattices).

Denote by $\tau (K,\gamma )$ the bottom of the essential spectrum of
the three-particle discrete Schr\"{o}dinger operator
$H_{\gamma}(K),\,K\in \T^{3}$ and by $N(K,\gamma;z)$ the number of
eigenvalues lying below $z\leq \tau (K,\gamma ).$

The main results of the present paper are as follows:

(i) for all $\gamma>0$ the operator $H_{\gamma}(0)$ has infinitely
many eigenvalues below the bottom of the essential spectrum and for
the number of eigenvalues $N(0,\gamma;z )$ lying below $z<0$ the
asymptotics
\[
\lim_{z\rightarrow -0}\frac{N(0,\gamma;z )}{|\log |z||}={\ U}
(\gamma ),\quad (0<{\ U} (\gamma )<\infty )
\]
holds, which is similar to the asimptotics founded in the continuous
case by Sobolev \eqref{sob};

(ii) for any $\gamma>0$ and $K\in \T^3\setminus \{0\}$
$N(K,\gamma;\tau (K,\gamma))$ is a finite number satisfying  the
following asymptotics
\[
\lim_{|K|\rightarrow 0}\frac{N(K,\gamma;0 )}{|\log |K||}=2{\ U}
(\gamma ).
\]
This result  is characteristic for the lattice system and does not
have any analogue in the continuous case.

We underline that these results are in contrast to similar results
for the continuous three particle Schr\"odinger operators, where the
number of eigenvalues does not depend on the three-particle total
momentum $K \in R^3.$

These results  are also in contrast with the results for
two-particle operators,in which  discrete Schr\"odinger operators
have finitely many eigenvalues for
 all $k \in U_\delta(0),$ where  $U_\delta(0)=\{K \in \T^3: |K|<\delta\}$
 is  a $\delta-$ neighborhood of the origin.

In the present paper  we provide an asymptotics
 for the Fredholm determinant $\Delta(k,0)$
resp. $\Delta(0,z)$ as $k\rightarrow0$ resp. $z\rightarrow0,$ which
plays a crucial role in the proof of the main results. In
particular, it has been  proved that $\Delta(k,0)$ resp.
$\Delta(0,z)$ is differentiable in $|k|$ at $k=0 \in T^3$ resp.in
$z$ at $z=0.$

The result concerning Efimov's effect has been proved in the
continuum case (see \cite{Sob,Tam94}) using resolvent expansion
established in \cite{JeKa}.

We note that  the operator $H_\gamma(K)$ has been considered before,
but only the existence of infinitely many eigenvalues below the
bottom of the essential spectrum of $H_\gamma(0)$ has been announced
in \cite{LaSh99} .

The organization of the present paper is as follows.

Section 1 is an introduction.

In Section 2 we introduce the Hamiltonians of systems of two and
three-particles in coordinate and momentum representations  as
bounded self-adjoint operators in the corresponding Hilbert spaces.

In Section 3 we introduce the total quasi-momentum and decompose the
energy operators into von Neumann direct integrals, choosing
relative coordinate systems.

In Section 4 we state the main results of the paper.

In Section 5 we study spectral properties of the two-particle
discrete Schr\"{o}dinger operator $\hbar_\gamma(k),k\in \T^3.$ We
prove the existence of a unique positive eigenvalue below the bottom
of the continuous spectrum of $\hbar_\gamma(k),k\in \T^3$ (Theorem
\ref{mavjud}) and obtain an asymptotics for Fredholm's determinant
associated with $\hbar_\gamma(k),k\in \T^3.$

In Section 6 we introduce the "channel operators" and describe its
spectrum by the spectrum of the two-particle discrete
Schr\"{o}dinger operators. Applying a Faddeev's  type integral
equation we establish the precise location of the essential spectrum
(Theorem \ref{ess}).

In Section 6 we  prove the finiteness of eigenvalues below the
bottom of the essential spectrum of $H(K)$ for $K\in \T^3\setminus
\{0\}$ (Theorem \ref{finite}).

In section 7 we  derive the asymptotics for the number of
eigenvalues $N(0,\gamma;z)$ resp. $N(K,\gamma;0)$ of $H_\gamma(0)$
as $z\rightarrow0$ resp. $H_\gamma(K)$ as $|K|\rightarrow 0$
(Theorem \ref{infinite}).

Throughout the present paper we adopt the following notations: We
denote by  $\T^3$   the three-dimensional torus, i.e., the cube
$(-\pi,\pi]^3$ with appropriately  identified sides. The torus
$\T^3$ will always be considered as an abelian group with respect to
the addition and multiplication by real numbers regarded as
operations on the three-dimensional space $\R^3$ modulo $(2\pi
\Z)^3$.

 For each  $\delta>0$ the notation
 $ U_\delta(0)=\{K \in \T^3: |K|<\delta\}$
 stands for a $\delta-$ neighborhood of the origin and
 $U^0_\delta(0)=U_\delta(0)\setminus \{0\} $ for a punctured
$\delta-$  neighborhood. The subscript $ \alpha $ (and also $\beta $
and $\gamma $) always runs from 1 to 3 and we use the convention  $
\alpha \not=\beta, \beta\not= \gamma , \gamma \not=\alpha  $.

\section{Description of the energy operators of two and
three particles on a lattice}

Let ${\Z}^{3}$ be the three-dimensional lattice and let $
\,\,(\Z^{3})^{m}\,,\,\,m\in N$ be the Cartesian $m-$ th power of $
\,\,\Z^{3}.\,\,$ Denote by $\ell _{2}((\Z^{3})^{3})$ the Hilbert
space of square-summable functions $\,\,\hat{\varphi}$ defined on
$(\Z^{3})^{3}$ and let $\ell _{2}^{a}((\Z^{3})^{3})\subset\ell
_{2}((\Z^{3})^{3})$ be the subspace  of functions antisymmetric with
respect to the first two coordinates.

We consider a system of three-particles (two identical fermions and
boson) moving on the three-dimensional lattice $\Z^3.$ Each fermions
are interact with the boson via a zero-range pair attractive
potential $\mu$. The free Hamiltonian $\,\,\widehat{H}^0_{\gamma
}\,\,$ of this system in the coordinate representation is an
operator on the Hilbert space $\,\,\ell _{2}^{a}((\Z^{3})^{3})\,\,$
of the form
\begin{align} \label{free}
&(\widehat{H}^0_{\gamma}\,\hat{\psi})(x_{1},x_{2},x_{3})=\no
\\
&\hspace{1cm}\sum_{s\in{\Z}^{3}}\hat{\varepsilon}(s)[\hat{\psi}(x_{1}+s,x_{2},x_{3})+
  \hat{\psi}(x_{1},x_{2}+s,x_{3})+
  \gamma \hat{\psi}(x_{1},x_{2},x_{3}+s)],\\
 & \hat{\psi} \in \ell
_{2}^{a}((\Z^{3})^{3})\no
\end{align}
where the function
 $\hat{\varepsilon}(s)$ defined on $\Z^{3}$ by
\begin{equation} \label{epsilon}
\hat{\varepsilon}(s)=\left\{
\begin{array}{ll}
3, & s=0 \\
-\frac{1}{2}, & \mid s\mid =1 \\
0, & \mbox{otherwise,}
\end{array}
\right.
\end{equation}
 the number $\gamma >0\,\,$ being the ratio of the mass of the
fermions to that of boson.

It is clear that the free Hamiltonian \eqref{free} is a bounded
self-adjoint operator on $\,\,\ell _{2}^{a}((\Z^{3})^{3}).$

The three-particle Hamiltonian $\widehat H_{\mu,\gamma}$ of the
quantum-mechanical three-particle
 system  is a bounded
perturbation of the free Hamiltonian $\widehat H^0_{\gamma}$
\begin{equation}\label{total}
 \widehat{H}_{\mu,\gamma}=\widehat{H}^0_{\gamma}-
\mu( \widehat{V}_{1}+\widehat{V}_{2}).
\end{equation}
Here $\widehat{V}_{\alpha}=\widehat{V}, \alpha=1,2 $ is a
multiplication operator on $\,\,\ell _{2}^{a}((\Z^{3})^{3})$
\begin{equation*}
(\widehat{V}_{\alpha}\hat{\psi})(x_1,x_2,x_3)= \delta_{x_\alpha
x_3}\hat{\psi}(x_1,x_2,x_3), \quad \hat{\psi} \in \ell
_{2}^{a}((\Z^{3})^{3}),
\end{equation*}
where $\delta_{x_\alpha x_3}$ is the Kroneker delta.

The Hamiltonian $\hat{h}_{\mu,\gamma }$
 describing a subsystem consisting  of a fermion  and
a boson is introduced as a bounded self-adjoint operator on the
Hilbert space $ \ell_{2}((\Z^{3})^{2})$
\begin{equation} \label{two.part}
\hat{h}_{\mu,\gamma}=\hat{h}^0_{\gamma}-\mu\hat{v},
\end{equation}
where
\begin{equation*}
(\hat{h}^0_{\gamma}\hat \psi)(x_\alpha,x_3)=
  \sum_{s\in{\Z}^{3}}\hat{\varepsilon}(s)[\hat{\psi}(x_{\alpha}+s,x_{3})+
  \gamma \hat{\psi}(x_{\alpha},x_{3}+s)],\quad \alpha=1,2,\,\hat{\psi}\in
  \ell _{2}((\Z^{3})^{2})
 \end{equation*}
 and
\begin{equation*}
(\hat{v}\hat{\psi})(x_\alpha,x_3)= \delta _{x_{\alpha}x_{3}}
\hat{\psi} (x_{\alpha},x_{3}),\quad \hat{\psi} \in
\ell_2(({\Z}^3)^2).
\end{equation*}

\subsection{ The momentum  representation}

Let $L^{a}_2((\mathbb{T}^{3})^{3})$ be the subspace of
$L_2((\mathbb{T}^{3})^{3})$ consisting of  functions which are
antisymmetric in the first two coordinates, where ${({\T}^3)^m}$
denotes the Cartesian $m$-th power of ${\T}^3=(-\pi,\pi]^3.$

Let ${ \cF}_m:L_2(({\T}^3)^m) \rightarrow \ell_2(( {\Z}^3)^m),\,m\in
\N$ denote the standard Fourier transform.

Denote by $\cF_3^{a}$ the restriction of $\cF_3$ to the subspace
$L^{a}_2((\T^{3})^{3})$. It is easy to check that $\cF_3^{a}:
L^{a}_2((\T^{3})^{3}) \rightarrow \ell^{a} _{2}((\Z^{3})^{2}).$ In
the momentum representation the three-resp. two-particle
Hamiltonians are given by the bounded self-adjoint operators on the
Hilbert space $L_2^{a}((\T^3)^3)$ resp. $L_2((\mathbb{T}^{3})^{2})$
\[
  H_{\mu,\gamma}= (\cF_3^{a})^{-1}\widehat{H}_{\mu,\gamma}{\cF_3^{a}}
\]
resp.
\[
  h_{\mu,\gamma}= (\cF_2)^{-1}\hat{h}_{\mu,\gamma}\cF_2.
\]

One has
\begin{equation} \label{Hamilt}
H_{\mu ,\gamma }=H^0_{\gamma }-\,\mu (V_{1}+V_{2}),
\end{equation}  where
\begin{equation}
 \label{free0}
 (H^0_{\gamma }\,f)(k_{1},k_{2},k_{3})=(\varepsilon
(k_{1})+\varepsilon (k_{2})+\gamma \,\varepsilon
(k_{3}))\,f(k_{1},k_{2},k_{3}),f\in L_2^{a}((\T^3)^3),
\end{equation}
\begin{align} \label{sob}
&(V_{\alpha
}f)(k_1,k_2,k_{3})=(V f)(k_1,k_2,k_{3}) \\
&=\frac{1}{(2\pi)^{3}}\int_{({\T}
^{3})^{3}}\delta(k_{\alpha}-k_{\alpha}^{\prime}) \delta(k_{\beta
}+k_{3}-k_{\beta}^{\prime}-k_{3}^{\prime})f(k_{1}^{\prime
},k_{2}^{\prime},k_{3}^{\prime})dk_{1}^{\prime }dk_{2}^{\prime
}dk_{3}^{\prime},\nonumber
\end{align}
$f\in L_{2}((\T^{3})^{3}),\quad \alpha ,\beta
=1,2,\quad \alpha \not=\beta $.\\
The function $\varepsilon $ is of the form
\begin{equation} \label{eps form}
\varepsilon (p)=\sum\limits_{i=1}^{3}(1-\cos {p^{(i)}}),\quad
p=(p^{(1)},p^{(2)},p^{(3)})\in {\T}^{3}
\end{equation}
and  $\delta (k)$ denotes the three-dimensional Dirac
delta-function.

The two-particle Hamiltonian $h_{\mu ,\gamma }$ is of the form:

\begin{equation} \label{two} h_{\mu ,\gamma }=h^0_{\gamma
}-\,\mu \,v,\\
\end{equation}
\[
  (h^0_{\gamma }\,f)(k_{\alpha},k_{3})=(\varepsilon (k_{\alpha})+
  \gamma
  \varepsilon (k_{3}))\,f(k_{\alpha},k_{3}), \quad  \alpha =1,2,\quad  f\in L_{2}((\T^{3})^{2})
\]
and the operator $v$ can be written in the form
\begin{align*}
  &(vf)(k_{\alpha},k_{3})=(v_\alpha f)(k_{\alpha},k_{3})\\
  &=\frac{1}{(2\pi)^{3}}\int_{({\T}^{3})^2}\delta
  (k_\alpha+k_3-k_\alpha^{\prime }-k_3^{\prime })f(k_\alpha^{\prime},k_3^{\prime})
  dk_\alpha^{\prime }dk_3^{\prime},\quad  \alpha =1,2,\quad
   f\in L_{2}(({\T}^{3})^2).
\end{align*}

\section{Decomposition of the energy operators into von Neumann direct integrals.
 Quasimomentum and coordinate systems}

 Given $m\in \N$, denote by $\hat U^m_s$, $s\in {\Z}^3,$
 the unitary operators on the Hilbert space
  $\ell_2(({\Z}^3)^m)$ defined by
\begin{equation*}
(\hat U^m_sf)(n_1,n_2,..., n_m)=f(n_1+s,n_2+s,...,n_m+s),\quad f\in
\ell_2(({\Z}^3)^m).
\end{equation*}

 We
easily see that
\begin{equation*}
  \hat U^m_{s+p}=\hat U^m_s\hat U^m_p,\quad s,p \in \Z^3,
  \end{equation*}
that is, $\hat U^m_s,s\in\Z^3 $ is a unitary  representation of the
abelian group $\Z^3$ in  $\ell_2(({\Z}^3)^3).$ Denote by $\hat
U^{m}_{as}$ the restriction of $\hat U^m_s,s\in\Z^3 $ to the
subspace $\ell ^{a}_2((\T^{3})^{3})$. Via the Fourier transform $
\cF_3^a$ the unitary representation $\hat U^{3}_{as}$ of $\Z^3$ in
$\ell_2^a(({\Z}^3)^3)$ induces the representation $ U^{3}_{as}$ of
the group $\Z^3$ in the Hilbert space $L_2^a(({\T}^3)^3)$  by the
unitary (multiplication) operators
\begin{align*}\label{grup}
&( U^{3}_{as} f)(k_1,k_2,k_3)=\exp \big (-i(s,k_1+k_2+k_3)\big
)f(k_1,k_2,k_3),\,\, s\in \Z^3,\,\,f\in L_2^a((\T^3)^3).
\end{align*}

For any $K\in {\T}^3$ we define $\F_K^3$ as follows
\begin{align*}
&\F_K^3=\{(k_1,k_{2},K-k_1-k_2){\in }({\T}^3)^3: \quad
k_1,k_2\in\T^3\}.\end{align*}

Let $L_2^{a} ( \F_K^3) $ be the subspace of $L_2( \F_K^3)$
consisting of functions which are antisymmetric in the first two
coordinates.

  Since the Hamiltonian  $ H_{\mu,\gamma}$ commutes with the group
   $U^{3}_{as}$, $s\in \Z^3$ the decomposition
\begin{equation*}
L^{(a)}_2((\T^3)^3)= \int_{K\in {\T}^3} \oplus L^{(a)}_2(\F_K^3)d K
\end{equation*}
of $ L^{(a)}_2((\T^3)^3)$  into the direct integral yields the
corresponding decompositions of the unitary representation
$U^{3}_{as}$,
 $s \in \Z^3$ and  the operator $H_{\mu,\gamma}$ into the  direct integral
\begin{align*}
&U^{3}_{as}= \int_{K\in {\T}^3} \oplus U^{3}_{as}(K)d K,\\
&H_{\mu,\gamma}=\int_ {K \in {\T}^3}\oplus \widetilde
H_{\mu,\gamma}(K)dK,
\end{align*}
where
\begin{equation*}
U^{(a)}_s(K)=\exp(-i(s,K))I \quad \text{on}\quad L^{(a)}_2(\F_K^3)
\end{equation*}
and $I=I_{L_2(\F_K^m)}$ denotes the identity operator on the Hilbert
space $ L^{(a)}_2(\F_K^m)$.

The Hamiltonian $ h_{\mu,\gamma}$ commutes with the group $ U^2_s$,
$s\in \Z^3$ and hence the operator $h_{\mu,\gamma}$ can be
decomposed into the direct integral
\begin{equation*}
 h_{\mu,\gamma}= \int_{k
\in {\T}^3}\oplus\tilde h_{\mu,\gamma}(k)d k,
\end{equation*}
with respect to  the decomposition
\begin{equation*}
 L_2 (( {\T}^3)^2) = \int_ {k \in
{\T}^3} {\ \oplus } L_2 ( \F_k^2 ) d k,
\end{equation*}
where
\begin{align*}
&\F_k^2=\{(k_1,k-k_{1}){\in }({\T}^3)^2: \quad
k_1\in\T^3\}.\end{align*} Now  we introduce the mappings
\begin{equation*}
\pi_{\alpha3}^{(3)}:(\T^3)^3\to (\T^3)^2,\quad
\pi_{\alpha3}^{(3)}((k_\alpha, k_\beta, k_3))=(k_\alpha, k_3),\quad
\alpha,\beta=1,2,\,\,\alpha \neq \beta,
\end{equation*}
and
\begin{equation*}
\pi_{\alpha}^{(2)}:(\T^3)^2\to \T^3,\quad
\pi_{\alpha}^{(2)}((k_\alpha, k_3))=k_\alpha.
\end{equation*}

Denote by $\pi^{(3)}_{K}$ , $K\in \T^3$ resp.
 $\pi^{(2)}_k$, $k\in
\T^3$ the restriction of $\pi^{(3)}_{\alpha 3}$ resp.
$\pi^{(2)}_{\alpha}$ onto $\F_K^3\subset (\T^3)^3$ resp.
$\F_k^2\subset (\T^3)^2$, that is,
\begin{equation}\label{project}
\pi^{(3)}_{K}=\pi^{(3)}_{\alpha3}\vert_{\F_K^3}\quad \text{and}\quad
\pi^{(2)}_{k}= \pi^{(2)}_{\alpha}\vert_{\F_k^2}.
\end{equation}
At this point it is useful to remark that $ \F^3_{K},\,\, K \in
{\bbT}^3 $ and $ \F^2_{k},\,\, k \in {\bbT}^3 $ are six and
three-dimensional manifolds isomorphic to ${({\bbT}^3)^2}$ and
${\bbT}^3,$ respectively.

\begin{lemma} The mappings $\pi_{K}^{(3)},\,K\in \T^3$
resp.$\pi_{k}^{(2)},\,k\in \T^3$ is bijective from
$\F_{K}^{3}\subset(\T^3)^3$ resp. $\F_{k}^{2}\subset(\T^3)^2$ onto
$(\T^3)^2$ resp. $\T^3$ with the inverse mappings given by
\[
  (\pi_{K}^{(3)})^{-1}(k_{\alpha},k_3)=
  (k_{\alpha},k_3,K-k_{\alpha}-k_3)
\]
resp.
\[
  (\pi_{k}^{(2)})^{-1}(k_{\alpha})=(k_{\alpha},k-k_{\alpha}).
\]
\end{lemma}
\qed

Let
\[
  U_{K}: L_{2}(\F_{K}^{3})\rightarrow L_{2}((\T^{3})^{2}),
  \quad U_{K}f=f\circ (\pi_{K}^{(3)})^{-1},\quad K\in \T^{3},
\]
and
\[
  u_{k}: L_{2}(\F_{k}^{2})\rightarrow L_{2}(\T^{3}),
  \quad u_{k}g=g\circ (\pi_{k}^{(2)})^{-1},\quad k\in \T^{3},
\]
where $\pi_{K}^{(3)}$ and $\pi_{k}^{(2)}$ are defined by
\eqref{project}. Then $U_{K}$ and $u_{k}$ are unitary operators.

Let
\[
  H_{\mu,\gamma}=U_{K}\widetilde{H}_{\mu,\gamma}(K)(U_{K})^{-1},
  \quad h_{\mu,\gamma}=u_{k}\tilde{h}_{\mu,\gamma}(k)(u_{k})^{-1}.
\]
The operator $h_{\mu,\gamma}(k),\, k\in \T^{3}$ is of the  form
\begin{align}\label{two-partSh}
  h_{\mu,\gamma}(k)=h^{0}_{\gamma}(k)-\mu v,
\end{align}
where
\begin{align}\label{two,multip}
&  h^{0}_{\gamma }(k)f(q)=\cE_{k,\gamma }(q)f(q),\qquad f\in L_{2}({\T}^{3}),\\
&\label{Ek}
 \cE_{k,\gamma}(q)=
  \varepsilon(q)+\gamma\varepsilon (k-q),
\end{align}
\begin{equation}\label{integ}
  (vf)(q)=
  \frac{1}{(2\pi)^{3}}\int_{{\T}^{3}}f(q^{\prime})dq^{\prime },
  \qquad f\in L_{2}({\T}^{3}).
\end{equation}

The operator $H_{\mu,\gamma}(K),\,K\in \T^{3}$ is of the form
\[
  H_{\mu,\gamma}(K)=H^{0}_{\gamma}(K)-2\mu V
\]
and in the coordinates $(k_\alpha,k_3)$ the operators
 $H^{0}_{\gamma
}(K)$ and $V$ are defined by
\begin{align*}\label{multip.oper}
  &(H^{0}_{\gamma
  }(K)f)(k_\alpha,k_3)=E(K,\gamma;k_\alpha,k_3)
  f(k_\alpha,k_3),
\quad f\in   L_{2}((\T^{3})^{2}),
\end{align*}
\begin{equation}\label{poten..}
  (Vf)(k_\alpha,k_3)=(V_\alpha f)(k_\alpha,k_3)=
  \frac{1}{(2\pi)^{3}}\int_{{\T}^{3}}f(k_\alpha,k_3^{\prime})
  dk_3^{\prime},\quad f\in L_{2}(({\T}^{3})^{2}),
 \end{equation}
where
\begin{equation}\label{EK}
  E(K,\gamma;k_\alpha,k_3)=\varepsilon(k_\alpha)+\varepsilon(k_3)+
  \gamma\varepsilon(K-k_\alpha-k_3),\quad k_\alpha,k_3 \in \T^3.
\end{equation}

\begin{remark} Since $v$ is positive its square root $v^{\frac{1}{2}}$
exists and the decomposition
$$L_{2}((\T^{3})^{2})= L_2(\T^{3})\otimes
L_2(\T^{3})$$ yields the representation $V^{\frac{1}{2}} =I\otimes
v^{\frac{1}{2}} $ for the operator $V^{\frac{1}{2}}.$
\end{remark}

\section{Statement of the main results}
In this section we give the precise formulation of the main results.

Set
\begin{equation}
  \mu_{0}=(2\pi)^{3}\,\big (\int_{{\T}^{3}}
  (\,\varepsilon(p)\,)^{-1}\,dp\big )^{-1}
\end{equation}
and fix uniquely
\begin{equation}\label{mu_0}
 v(\gamma)=\mu_{0}(1+\gamma).
 \end{equation}

 We introduce the
family of the two-particle operators $\hbar_{\gamma}(k)$ depending
on $k\in{\T}^{3},\gamma>0$
\begin{equation} \label{hlambda}
 \hbar_{\gamma}(k)\equiv h_{v(\gamma),\gamma}(k).
\end{equation}

\begin{definition} \label{resonance} The operator $\hbar_{\gamma}(0)$
is said to have a zero energy resonance if the number $1$ is an
eigenvalue for the operator
\begin{equation}\label{equation}
  {G_{\gamma}}=v^{\frac{1}{2}}(\hbar_{\gamma}(0))^{-1}v^{\frac{1}{2}}
  \end{equation}
and the associated eigenfunction  $\psi$ satisfies the following
condition $(v^{1/2}\psi)(0)\neq 0.$ Without loss of generality we
can always normalize $(v^{1/2}\psi)(0)$ so that
$(v^{1/2}\psi)(0)=1.$
\end{definition}

\begin{remark}\label{resonarem}
If $\psi $
 is eigenfunction of $G_{\gamma}$ associated to eigenvalue $1$
 then  the
function
$$
 f (p)={(v^{1/2}\psi)(p)}{(\varepsilon (p))^{-1}}
$$
is a simple solution (up to a constant factor) of the equation
$\hbar_{\gamma }(0)f=0.$

If $v^{1/2}\psi(0)\neq0$ from
$$\varepsilon(p)=\frac{1}{2}|p|^{2}+O(|p|^{4})\quad \mbox{ as}\quad
p\rightarrow 0$$
 we conclude that the operator $\hbar_{\gamma }(0)$ has a zero
 energy resonance precisely if $G_\gamma$ has 1 as eigenvalue and the
 associated eigenfunction $\psi$ satisfy $v^{1/2}\psi(0)\neq 0.$
\end{remark}
Let $H_{\gamma}(K)$ be the family of the three-particle operators
defined by
  \begin{equation} \label{Hlambdaa}
  H_{\gamma}(K)\equiv H_{v(\gamma),\gamma}(K).
\end{equation}

We recall that the main goal of the paper is to study spectral
properties
 of the operators   $ H_{\gamma}(K) $ depending on the parameters $K\in
\T^{3}$ (the three-particle quasi-momentum) and $\gamma>0$ (ratio of
the mass of the fermions and the boson).

  For each
$K\in \T^{3}$ and $\gamma>0$ we set:
\[
  E_{\min }(K,\gamma )=\min_{p,q \in {\T}^{3}} E(K,\gamma;p,q),\quad
  E_{\max}(K,\gamma )=\max_{p,q \in {\T}^{3}} E(K,\gamma;p,q).
\]

The main results of the present paper are as follows.

\begin{theorem}\label{mavjud}
 For any  $\gamma >0$ the operator
 $\hbar_{\gamma}(0)$ has a zero-energy resonance.
  For all $k\in\T^{3}\setminus \{0\}$ the operator $\,\hbar_{\gamma
}(k)\,$ has a unique positive eigenvalue $z_{\gamma}(k)$ below the
bottom of the essential spectrum of $\,\hbar_{\gamma }(k)$ and
$z_{\gamma}(k)$ is even and real-analytic in $\T^{3}\setminus
\{0\}.$
\end{theorem}

Let $\tau (K,\gamma )=\inf_{p \in {\T}^{3}}\{z_{\gamma}(K-p)+
\varepsilon(p)\}.$  We denote by $N(K,\gamma;z )\,$ the number of
eigenvalues of $\,H_{\gamma }(K)\,$ below $\,z\leq \tau (K,\gamma
).$

\begin{theorem}\label{ess} For the essential spectrum
$\sigma _{ess}(H_{\gamma }(K))$ of $\,\,H_{\gamma }(K)$ the
following equality
\begin{equation}\label{esss}
  \sigma_{ess}(H_{\gamma}(K))=[\tau(K,\gamma)
  ,E _{max}(K,\gamma )]
\end{equation}
holds.
\end{theorem}

In the following theorems we   describe precisely the dependence of
the number of eigenvalues of $ H_{\gamma}(K) $ lying below the
bottom of the essential spectrum  $\tau(K,\gamma) $ on the
parameters $K\in \T^{3}$ and $\gamma>0$.

\begin{theorem}\label{finite}
 For any  $\gamma
>0$ and   $K\in \T^3\setminus\{0\}$ the operator $H_{\gamma }(K)$ has
a finite number of eigenvalues lying below  $\tau(K,\gamma) $ .
\end{theorem}
\begin{theorem}\label{infinite}
For any $\gamma >0$ the operator $H_{\gamma}(0)$ has infinitely many
eigenvalues lying below the bottom of the essential spectrum. The
function $N(0,\gamma;z )$ resp. $N(K,\gamma;0 )$ obeys the relation
\begin{equation}\label{asympz}
  \lim_{z\rightarrow 0-}\frac{N(0,\gamma;z )}
  {\mid \log \mid z\mid \mid }\,=
  {\ U}(\gamma )>0,
\end{equation}
resp.
\begin{equation}\label{asympK}
  \lim_{\mid K\mid \rightarrow 0  }\frac{N(K,\gamma;0)}{\mid \log
  \mid K\mid \mid }\,=2\,{\ U}(\gamma ).
\end{equation}
\end{theorem}
\begin{remark}
Clearly, the infinitude of the negative discrete spectrum of
$H_{\gamma}(0)$  follows automatically from the strict positivity of
${\ U}(\gamma ).$
\end{remark}

\section {Spectral properties of the two-particle operator
              $ \hbar_{\gamma }(k)$}

In this section we study some spectral properties of the
two-particle discrete Schr\"{o}dinger operator $ \hbar_{\gamma
}(k),\, k\in {\T}^{3}$ defined by \eqref{hlambda}.

The perturbation $v$ of the multiplication operator $h^0_\gamma(k)$
is a bounded self-adjoint operator of rank one. Therefore in
accordance to Weyl's theorem  the essential spectrum of
$\hbar_{\gamma }(k), \, k\in {\T}^{3},$ fills the following interval
on the real axis:
\[
  \sigma_{ess}(\hbar_{\gamma }(k))=[\,m(k,\gamma)\,,\,M(k,\gamma )\,],
\]
where
\[
  m(k,\gamma )=\,\min_{p\in {\T}^{3}}\cE_{k,\gamma}(p),\qquad
  M(k,\gamma )=\,\max_{p\in {\T}^{3}}\cE_{k,\gamma }(p)
\]
and the function $\cE_{k,\gamma }(p)$ is defined by \eqref{Ek}.

 Let $\C$ be the field of complex numbers. For any $k\in \T^{3}$ and $z\in
\C \setminus [\,m(k,\gamma)\,,\,M(k,\gamma )\,],$ we define a
function (the Fredholm determinant associated with the operator
$\hbar_{\gamma }(k)$)
\begin{equation}\label{det}
  \Delta_{\gamma }(k,z)= 1-\frac{v(\gamma)}
  {(2\pi)^{3}}\int_{\T^{3}}
  (\cE_{k,\gamma}(q)-z)^{-1}dq,
\end{equation}
where $v(\gamma)$ is defined in \eqref{mu_0} and $\gamma>0$.

Note that the function $\Delta_{\gamma }(k,z)$ is real-analytic in
$\T^{3}\times(\C\backslash [m(k,\gamma );M(k,\gamma )]). $

The following lemma is a simple consequence of the Birman-Schwinger
principle and  Fredholm's theorem.

\begin{lemma}\label{nollar} For any $k\in \T^{3}$ and $\gamma>0$ a point
$z\in \C\setminus [\,m(k,\gamma)\,,\,M(k,\gamma )\,],$ is an
eigenvalue of the operator $\hbar_{\gamma}(k)$ if and only if $
  \Delta_{\gamma}(k,z)=0.$
\end{lemma}$\Box$

 Let
\begin{equation}\label{det}
  \Delta_{\gamma }(0,0)= 1-\frac{v(\gamma)}
  {(2\pi)^{3}}\int_{\T^{3}}
  \cE_{0,\gamma}(q))^{-1}dq.
\end{equation}

Note that
$$
\lim_{z\to 0-} \Delta_{\gamma }(0,z)=\Delta_{\gamma }(0,0).
$$
The following Lemma is evident
\begin{lemma}\label{resonance}(see\cite{ALzM04}). For any $\gamma >0$ the operator
$\hbar_{\gamma}(0)$  has a zero-energy resonance if and only if
$\Delta_{\gamma}(0,0)=0.$
\end{lemma}$\Box$

{\it {Proof of Theorem \ref{mavjud}.}}

for any $\gamma>0$ for the equation
\begin{equation}\label{equation}
  {G_{\gamma}}=v^{\frac{1}{2}}(\hbar_{\gamma}(0))^{-1}v^{\frac{1}{2}}
  \end{equation}
the $const.\neq0$ is the nonzero solution and hence by Definition
\ref{resonance}
 we
can conclude that the operator
 $\hbar_{\gamma}(0)$ has a zero-energy resonance.

Case 1. Let either $\gamma\neq 1$ and $k\in \T^3$ or $\gamma=1$ and
$k\in (-\pi,\pi)^3$.

 The function $\cE_{k,\gamma}(p)$ can
be rewritten in the form
\begin{equation}\label{Eraz2}
  \cE_{k,\gamma}(p)=3(1+\gamma)-\sum\limits_{j=1}^{3}\sqrt{1+2\gamma\cos{k^{(j)}}
  +\gamma^{2}}\,\cos(p_{j}-p_{\gamma}(k^{(j)})),
\end{equation}
where \begin{equation}\label{pk}
  p_{\gamma}(k^{(j)})=\arcsin\frac{\gamma\sin{k^{(j)}}}{\sqrt{1+2\gamma\cos{k^{(j)}}
  +\gamma^{2}}},\quad k^{(j)}\in (-\pi,\pi],\quad j=1,2,3.
\end{equation}
Taking into account \eqref{Eraz2} we have that the vector-function
\[
  p_{\gamma}:\T^{3}\rightarrow \T^{3},\quad
  p_{\gamma}(k)=p_{\gamma}(k^{(1)},k^{(2)},k^{(3)})=
  (p_{\gamma}(k^{(1)}),p_{\gamma}(k^{(2)}),p_{\gamma}(k^{(3)}))\in\T^{3}
\]
is odd and  regular in $(-\pi,\pi)^3$ and $\cE_{k,\gamma}(p)$
reaches the minima on it. One has, as easily seen from the
definition,
\begin{equation}\label{pkasymp}
  p_{\gamma}(k)=\frac{\gamma}{1+\gamma}k+O(|k|^{3})\quad
  \mbox{as}
  \quad k\rightarrow 0.
\end{equation}

 Moreover from \eqref{Eraz2}
it follows that
\begin{equation}\label{minE}
  m(k,\gamma)=\cE_{k,\gamma}(p_{\gamma}(k))=
  3(1+\gamma)-\sum\limits_{j=1}^{3}\sqrt{1+2\gamma\cos{k^{(j)}}+\gamma^{2}}.
\end{equation}

By Lemma \ref{resonance}
\[
  \Delta_{\gamma}(0,0)=1-\frac{\mu_{0}}{(2\pi)^{3}}
  \int_{{\T}^{3}}(\varepsilon(q))^{-1}dq=0.
\]
Since for any $\gamma>0$ the function $\Delta_{\gamma}(0,z)$ is
monotone decreasing on $(-\infty,m(k,\gamma))$ for any $z<0$ the
inequality $\Delta_{\gamma}(0,z)>0$ holds. By Lemma \ref{nollar} the
operator $\hbar_{\gamma}(0)$ does not have any  negative spectrum.
Thus the operator $\hbar_{\gamma}(0)$ is positive.

 Since
$p=p_{\gamma}(k)$ is the non-degenerate minimum of the function
$\cE_{k,\gamma}(p)$ the function
$(\cE_{k,\gamma}(p)-m(k,\gamma))^{-1}$ is integrable and we define
$\Delta_{\gamma}(k,m(k,\gamma))$ by
\[
  \Delta_{\gamma}(k,m(k,\gamma))
  =1-\frac{(1+\gamma)\mu_{0}}{(2\pi)^{3}}
  \int_{{\T}^{3}}(\cE_{k,\gamma}(q)-m(k,\gamma))^{-1}dq.
\]
By dominated convergence theorem we have
\[
  \lim_{z\rightarrow m(k,\gamma)-0} \Delta_{\gamma}(k,z)=
  \Delta_{\gamma}(k,m(k,\gamma)).
\]
From the representations \eqref{Eraz2} and \eqref{minE} it follows
that for all $k\neq 0,\, q\neq 0$ the inequality
\[
  \cE_{k,\gamma}(q+p_{\gamma}(k))-m(k,\gamma)
  <\cE_{0,\gamma}(q)
\]
holds and hence we obtain the following inequality
\begin{equation}\label{detk<0}
 \Delta_{\gamma}(k,m(k,\gamma))
 <\Delta_{\gamma}(0,0)=0,\quad k\neq 0.
\end{equation}
For each $\gamma>0$ and $k\in{\T}^{3}$ the function
$\Delta_{\gamma}(k,\cdot)$ is continuous monotone decreasing on
$(-\infty,m(k,\gamma))$ and $\Delta_{\gamma}(k,z)\rightarrow 1$ as
$z\rightarrow -\infty.$ Therefore by virtue of  \eqref{detk<0} there
is a  unique number $z_{\gamma}(k)\in(-\infty,m(k,\gamma))$ such
that $\Delta_{\gamma}(k,z_{\gamma}(k))=0.$  By Lemma \ref{nollar}
for any nonzero $k\in {\T}^{3}$ the operator $\hbar_{\gamma}(k)$ has
a unique eigenvalue below $m(k,\gamma).$ For any $\gamma>0$ and
$z\in(-\infty,m(k,\gamma)]$ the equality
$\Delta_{\gamma}(-k,z)=\Delta_{\gamma}(k,z),$ $k\in \T^3$
 holds.
Hence $z_{\gamma}(k)$ is even.

Let us prove the positivity of the eigenvalue $z_{\gamma}(k),\,
k\neq 0.$ First we verify, for all $k\in\T^{3},\quad k\neq 0,$ the
inequality
\begin{equation}\label{5.11}
  \Delta_{\gamma}(k,0)>0.
\end{equation}

We have
\begin{equation}
  \Delta_{\gamma}(k,0)=\frac{\mu_0}{(2\pi)^{3}}
  \int_{\T^{3}}
  \frac{\gamma (\varepsilon(k-q)-\varepsilon(q))}{\varepsilon(q)
 \cE_{k,\gamma}(q)}dq.
\end{equation}
Making a change of variables $q=\frac{k}{2}-p$ in (5.12) and using
the equality $\Delta_{\gamma}(k,0)=\Delta_{\gamma}(-k,0)$ it is easy
to show that
\[
  \Delta_{\gamma}(k,0)=\frac{\Delta_{\gamma}(k,0)+\Delta_{\gamma}(-k,0)}{2}=
\]
\[
  =\frac{\mu_0}{2(2\pi)^{3}}\int_{\T^{3}} \gamma
  (\varepsilon(\frac{k}{2}+p)-\varepsilon(\frac{k}{2}-p))^{2}\,F(k,p)dp,
\]
where
\[
  F(k,p)=\frac{\varepsilon(\frac{k}{2}+p)+\varepsilon(\frac{k}{2}-p)}
  {\varepsilon(\frac{k}{2}+p)\varepsilon(\frac{k}{2}-p)
  \cE_{k,\gamma}(\frac{k}{2}+p)\cE_{k,\gamma}(\frac{k}{2}-p)}>0.
\]
Thus the inequality \eqref{5.11} is proven.

For any $\gamma>0$ and $k\in\T^{3}$ the function
$\Delta_{\gamma}(k,\cdot)$ is monotone decreasing and the
inequalities
\[
  \Delta_{\gamma}(k,0)>\Delta_{\gamma}(k,z_{\gamma}(k))=0>
  \Delta_{\gamma}(k,m(k,\gamma)),\quad k\neq 0
\]
hold. Therefore the eigenvalue $z_{\gamma}(k)$ of the operator
$\hbar_{\gamma}(k)$ belongs to the interval $(0,m(k,\gamma)).$

Case 2. Let $\gamma=1$ and $k\in \T^3\setminus (-\pi,\pi)^3$. In
this case if we take into account the equality
$$\frac{(1+\gamma)\mu_{0}}{(2\pi)^{3}}
  \int_{{\T}^{3}}\cE_{k,\gamma}(q)-m(k,\gamma))^{-1}dq=-\infty,
  $$
then the proof of  Theorem \ref{mavjud} going on by the same way as
 Case 1.

The function $\Delta_{\gamma}(k,z)$ is analytic in $(k,z)\in
(\T^3)\times (-\infty,0) $ and hence $z(k)$ as a unique solution of
the equation $$\Delta_{\gamma}(k,z)=0$$ is analytic in $\T^3.$

 $\Box$

The following decomposition is important for the proof of the main
results \eqref{asympz} and \eqref{asympK}.
\begin{lemma}\label{detraz}
For any $\gamma>0$ and $k\in \T^3$  and $z\leq m(k,\gamma)$ the
following decomposition holds:
\[
  \Delta_{\gamma}(k,z)=\frac{v(\gamma)\,
  \gamma^{3/2}}{\sqrt{2}\,\pi(1+\gamma)^{3/2}}
  \left[m(k,\gamma)-z\right]^{\frac{1}{2}}+\Delta_{\gamma}^{(20)}
  (m(k,\gamma)-z),
\]
where $\left[m(k,\gamma)-z\right]^{\frac{1}{2}}>0$ for
$m(k,\gamma)-z>0$ and
$\Delta_{\gamma}^{(20)}(m(k,\gamma)-z)=O(m(k,\gamma)-z)$ as
$z\rightarrow m(k,\gamma).$
\end{lemma}
\begin{proof}
Let
\[
  E_{\gamma}(k,p)=\cE_{k,\gamma}(p+p_{\gamma}(k))-m(k,\gamma).
\]
 Then using \eqref{Eraz2} we
conclude
\[
  E_{\gamma}(k,p)=\sum_{j=1}^{3}\sqrt{1+2\gamma \cos
  k_{j}+\gamma^{2}}\,(1-\cos p_{j}).
\]

We define the function $\widetilde{\Delta}_{\gamma}(k,\omega)$ on
$\T^{3}\times \C_{+}$ by
$\widetilde{\Delta}_{\gamma}(k,\omega)=\Delta_{\gamma}
(k,m(k,\gamma)-\omega^2),$ where $\C_{+}=\{z\in \C:\quad \Re z>0\}.$
The function $\widetilde{\Delta}_{\gamma}(k,\omega)$ can be
represented in the following way:
\[
  \widetilde{\Delta}_{\gamma}(k,\omega)=1-v(\gamma)(2\pi)^{-3}
  \int_{\T^{3}}\,\frac{dp}{E_{\gamma}(k,p)+\omega^{2}}=
\]
\[
  =1-v(\gamma)(2\pi)^{-3}\int_{\T^{3}}\,\frac{dp}
  {\sum_{j=1}^{3}\sqrt{1+2\gamma \cos
  k_{j}+\gamma^{2}}\,(1-\cos p_{j})+\omega^{2}}.
\]

Let $V_{\delta}(0)$ be the complex $\delta -$ neighborhood of the
point $\omega = 0 \in\C.$ Denote by
${\Delta}_{\gamma}^{\ast}(k,\omega)$ the analytic continuation of
the function $\widetilde{\Delta}_{\gamma}(k,\omega)$ to the region
$\T^{3}\times (\C_{+}\cup V_{\delta}(0))$(see \cite{Ltmf92}). The
function $\widetilde{\Delta}_{\gamma}(\cdot,\omega)$ is even in
$k\in \T^{3}.$ A Taylor series expansion gives
\[
  {\Delta}_{\gamma}^{\ast}(k,\omega)=
  \widetilde{\Delta}_{\gamma}^{(01)}(k,0)\,\omega+
  \widetilde{\Delta}_{\gamma}^{(02)}(k,\omega)\,\omega^{2},
\]
where $\widetilde{\Delta}_{\gamma}^{(02)}(k,\omega)=O(1)$ as
$\omega\rightarrow 0.$ Then a simple computation shows that
\begin{equation}\label{dere}
  \frac{\partial{\Delta}_{\gamma}^{\ast}(0,0)}{\partial\omega}=
  \widetilde{\Delta}_{\gamma}^{(01)}(0,0)=
  \frac{v(\gamma)\,\gamma^{3/2}}{\sqrt{2}\,\pi(1+\gamma)^{3/2}}\neq 0.
\end{equation}
 The equality $m(k,\gamma)-z=\omega^2$ yields the proof of Lemma.
\end{proof}
\begin{corollary}\label{corol}
The function $z_\alpha(k)=m(k,\gamma)-w^2_\gamma(k)$ is
real-analytic in $\T^3,$ where $w_\alpha(k)$ is a unique simple
solution of the equation $ {\Delta}_{\gamma}^{\ast}(k,\omega)=0$ and
$w_\alpha(k)=O(|k|^2)$ \text{as}\quad $k\to 0$.
\end{corollary}
\begin{proof}
The equation $\tilde\Delta_\alpha(k,w)=0$ has a unique simple
solution $w_\alpha(k),\,k\in \T^3$ and it is real-analytic in $\in
\T^3.$ Taking into account that the function
$\tilde\Delta_\gamma(k,w)$ is even in $k\in U_\delta(0),\delta>0$
and $w_\gamma(0)=0$ we have that $w_\gamma(k)=O(|k|^2).$ Therefore
the function $z_\gamma(k)=E^{(\gamma)}_{\min}(k)-w^2_\gamma(k)$ is
real-analytic in $U_\delta(0).$
\end{proof}
 \begin{lemma}\label{tasvir}
 For any $k \in \T^3\setminus\{0\}$ there exists a number
 $\delta(k)>0$ such that, for all $z \in
 V_{\delta(k)}(z_\gamma(k))$, where $V_{\delta(k)}(z_\gamma(k))$
 is the  $\delta(k)$-neighborhood of the point  $z_\gamma(k)$,
 the following representation holds
$$
\tilde\Delta_\gamma(k,z)=C_1(k)(z-z_\gamma(k))\hat\Delta_\gamma(k,z).
$$
Here $C_1(k)\neq 0$ and $\hat\Delta_\gamma(k,z)$ is regular in
$V_{\delta(k)}(z_\gamma(k))$ and
$\hat\Delta_\gamma(k,z_\gamma(k))\neq 0.$
\end{lemma}

\begin{proof}
Since $z_\gamma(k)<m(k,\gamma),\,k\neq0$ the function
$\Delta_{\gamma}(k,z)$ is regular in $\C\backslash [m(k,\gamma
);M(k,\gamma )].$ Hence for some $\delta(k)>0$ it
 can be expanded as
 $$
\Delta_\gamma(k,z)=\sum_{n=1}^{\infty}C_n(k)(z-z_\gamma(k))^n, \quad
z\in V_{\delta(k)}(z_{\gamma}(k)),$$ where
\begin{align*}
&C_1(k) = \frac{\mu_\gamma^o m_{\beta\gamma}^{3/2}}{ \sqrt{2}\pi}
\frac{1}{2\sqrt{m(k,\gamma)-z_\gamma(k)}}\neq 0, \quad k\neq 0.
\end{align*}

Clearly, that $\hat\Delta_\gamma(k,z)$ is regular in
$V_{\delta(k)}(z_{\gamma}(k))$.
 Since $z_\gamma(k),\, k\neq 0$ is a unique simple solution of the
equation $\Delta_\gamma(k,z)=0,\,z< m(k,\gamma)$, we have
$\hat\Delta_\gamma(k,z_\alpha(k))\neq 0.$
\end{proof}

\section{ The essential spectrum  of
 \ $H_{\gamma}(K)$ and channel operators}

Recall that we consider a three-particle system  consisting of two
identical fermions and boson.  The fermions interact with boson via
a zero-range pair attractive potential. Therefore we have only one
non-trivial channel operator $H_{\gamma }^{ch}(K),\,K\in \T^{3}$
acting on $ L_{2}((\T^{3})^{2})$ as
\begin{equation}
  (H_{\gamma}^{ch}(K)f)(p,q)=E(K,\gamma;p,q)f(p,q)-
  \frac{v(\gamma)}{(2\pi)^{3}}\int_{{\T}^{3}}f(p,q^{\prime})dq^{\prime},
\end{equation}
where  $E(K,\gamma;p,q)$ is defined by \eqref{EK}.

Since the operator $H_{\gamma }^{ch}(K)$ commutes with the group
$\{U_{s}^{(2)},\, s\in \Z^{3}\}$ of the unitary operators
\[
  (U_{s}^{(2)}f)(p,q)=exp\{-\mathrm{i}(s,p)\}\,f(p,q),\quad f\in L^{(a)}_{2}((\T^{3})^{2})
\]
the decomposition of the Hilbert space $L_{2}((\T^{3})^{2})$ into
the direct integral
\begin{equation*}
  L_{2}(({\T}^{3})^{2})= \int_{{\T}^{3}}\,\oplus \,L_{2}({\T}^{3})\,dp
\end{equation*}
yields the decomposition
\begin{equation*}
   H_{\gamma}^{ch}(K)=\int_{{\T}^{3}}\oplus\,
   H_{\gamma}^{ch}(K,p)\,dp.
\end{equation*}

The fiber operator $H_{\gamma}^{ch}(K,p)$ acts in  $L_{2}(\T^{3})$
by
\begin{equation}\label{channel}
  H_{\gamma}^{ch}(K,p)=\hbar_{\gamma}(K-p)+
  \varepsilon(p)\,I_{L_2(\T^3)}
\end{equation}
where $\,I_{L_2(\T^3)}\,$ is identity operator on $L_2(\T^3)$ and $
\hbar_{\gamma }(k) $ is the two-particle operator defined by
\eqref{hlambda}.

The representation \ref{channel} of the operator
$\,H_{\gamma}^{ch}(K,p)\,$ implies the equality
\begin{align}\label{6.5}
  &\sigma(H_{\gamma}^{ch}(K,p))=\{z_{\gamma}(K-p)+ \varepsilon(p)\} \cup\\
 &\hspace{4cm}  [ m(K-p,\gamma )+
  \varepsilon(p),M(K-p,\gamma )+
  \varepsilon(p)].\nonumber
\end{align}
\begin{remark}
We note that the point $ z_{\gamma}(0)+ \varepsilon(0)=0 $ is a zero
energy resonance of the operator $H_{\gamma}^{ch}(0,0).$
\end{remark}

The Theorem (see, e.g., \cite{RSIV}) on the spectrum of decomposable
operators and the structure \eqref{6.5} of the spectrum of
$H_{\gamma}^{ch}(K,p)$ obtained above lead to the following results:

\begin{lemma}\label{chan,cpec} The equality
\begin{equation}\label{ch.sp.str}
  \sigma(H_{\gamma}^{ch}(K))=
[\tau(K,\gamma),E_{max}(K,\gamma)]
\end{equation}
holds, where  $E_{max}(K,\gamma)$ is the maximum value of the
function $E(K,\gamma;p,q)$.
\end{lemma}

\begin{lemma}\label{tau<E} For any $K\in T^3$ the following inequality
\[
  {\tau}(K,\gamma )<E_{\min}(K,\gamma)
\]
holds.
\end{lemma}

\begin{theorem}\label{ChanEss}
For the essential spectrum $\sigma_{ess}(H_{\gamma }(K))$ of
$H_{\gamma }(K)$ the equality
\[
  \sigma(H_{\gamma}^{ch}(K))=\sigma_{ess}(H_{\gamma }(K))
\]
holds.
\end{theorem}
The proof of Theorem \ref{ChanEss} is similar to the one of Theorem
4.3
 in \cite{ALzM04}.

{\bf Proof of Theorem \ref{ess}.}\marginpar{yozish kerak} Theorem
\ref{ess} follows easily from Lemma \ref{chan,cpec} and Theorem
\ref{ChanEss}. \qed

Let $W_\gamma(K,z),K\in \T^{3},\,z < {\tau}(K,\gamma )$  be the
operators on $L_2(({\T}^3)^2)$ defined as
$$
W_\gamma(K,z)=I+V^{\frac{1}{2}} R^{ch}_\gamma (K,z)V^{\frac{1}{2}},
 $$
where  $R^{ch}_\gamma (K,z)$ is the resolvent of
$H_{\gamma}^{ch}(K).$ One checks that
$$
W_\gamma(K,z)=(I-V^{\frac{1}{2}}R^{0}_\gamma (K,z)V^{\frac{1}{2}})^{-1},$$
where $R^{0}_\gamma (K,z)$ the resolvent of the operator $H_0^\gamma(K).$

Note that for all $K\in \T^{3},\,z < {\tau}(K,\gamma )$ the
operators $W_\gamma(K,z)$ are positive.

Denote by $$
 {\bf T_\gamma}(K,z),\,K\in \T^{3},\,  z < {\tau}(K,\gamma )
$$
the operator in $ L^{(a)}_2(({\T}^3)^2)$ defined by

\begin{equation}
 {\bf T} _\gamma( K, z) = 2W^{\frac{1}{2}}_\gamma(K,z)V^{\frac{1}{2}}
R^{0}_\gamma (K,z)V^{\frac{1}{2}}W^{\frac{1}{2}}_\gamma(K,z).
\end{equation}
For any bounded self-adjoint operator $A$ acting in the Hilbert
space ${\cH}$ not having any essential spectrum on the right of the
point $z$ we denote by  ${\cH}_A(z)$ the subspace such that $(Af,f)
> z(f,f)$ for any $f \in {\cH}_A(z)$ and set
$n(z,A)=\sup_{\cH_A(z)}\dim{\cH}_A(z)$.

By the definition of $N(K, \gamma,z)$ we have
$$
N(K,\gamma,z)=n(-z,-H_\gamma(K)),\,-z > -{\tau}(K,\gamma ).
$$
The following lemma is a realization of the well known
Birman-Schwinger principle for the three-particle Schr\"{o}dinger
operators on a lattice (see \cite{Sob,Tam94} ).
\begin{lemma}\label{b-s}
 For $z<{\tau}(K,\gamma )$ the operator ${\bf T} _\gamma(K,z)$ is
compact and continuous in $z$ and
$$
N(K, \gamma, z)=n(1,{\bf T}(K,z)).
$$
\end{lemma}
\begin{proof} We first verify the equality
\begin{equation}\label{tenglik}
N(K,\gamma,z)=n(1,2(R^0_\gamma(K,z))^{\frac{1}{2}}V(R^0_\gamma(K,z))^{\frac{1}{2}}).
\end{equation}
 Assume that $u \in
{\cH}_{-H_\gamma(K)}(-z)$, that is, $((H^0 _\gamma(K)-z)u,u) < 2( Vu,u).$
Then
$$ (y,y) < 2(R^{\frac{1}{2}}_0(K,z) V
(R^0_\gamma(K,z))^{\frac{1}{2}}y,y),\quad y=(H^0 _\gamma(K)-z)^{\frac{1}{2}}u.
$$
Thus $N(K, \gamma,z) \leq
n(1,2(R^0_\gamma(K,z))^{\frac{1}{2}}V(R^0_\gamma(K,z))^{\frac{1}{2}})$. Reversing the
argument we get the opposite inequality, which proves
\eqref{tenglik}. Any nonzero  eigenvalue of
$(R^0_\gamma(K,z))^{\frac{1}{2}}V^{\frac{1}{2}}$ is an eigenvalue for
$V^{\frac{1}{2}}(R^0_\gamma(K,z))^{\frac{1}{2}}$
 as well, of the same
algebraic and geometric multiplicities.

Therefore  we get
$$
n(1,2R^{\frac{1}{2}}_0(K,z) V R^{\frac{1}{2}}_0(K,z))=
n(1,2V^{\frac{1}{2}} R_0(K,z)V^{\frac{1}{2}}).
$$

Let us check that
$$
n(1,2(R^0_\gamma(K,z))^{\frac{1}{2}} V(R^0_\gamma(K,z))^{\frac{1}{2}})=n(1,{\bf
T_\gamma}(K,z)).
$$
We shall show that for any $u \in {\cH}_{2(R^0_\gamma(K,z))^{\frac{1}{2}}V
(R^0_\gamma(K,z))^{\frac{1}{2}}}(1)$ there exists $y \in{\cH}_{{\bf
T} _\gamma (K,z)}(1)$ such that $(y,y)<({\bf T} _\gamma (K,z)y,y).$ Let $u \in
{\cH}_{2(R^0_\gamma(K,z))^{\frac{1}{2}}V (R^0_\gamma(K,z))^{\frac{1}{2}}}(1)$ that is,
$$
(u,u)< 2(V^{\frac{1} {2}} R^0_\gamma(K,z) V^{\frac{1}{2}} u,u)
$$
and hence
\begin{equation}\label{coordinate}
((I-V^{\frac{1}{2}} R^0_\gamma(K,z) V^{\frac{1}{2} })u,u)<(V^{\frac{1}{2}
}R^0_\gamma(K,z)V^{\frac{1}{2}} u,u).
\end{equation}
Denoting by $y=(I- V^{\frac{1}{2}}R^0_\gamma(K,z)V
^{\frac{1}{2}})^{\frac{1}{2}} u $ we have
$$
(y,y)< 2(W^{\frac{1}{2}} _\gamma(K,z)V^{\frac{1}{2}} R^0_\gamma(K,z)
V^{\frac{1}{2}} W^{\frac{1}{2}}  _\gamma(K,z)y,y),
$$
that is, $ (y,y)\leq ({\bf T} _\gamma(K,z)y,y). $ Thus $
n(1,2(R^0_\gamma(K,z))^{\frac{1}{2}}V (R^0_\gamma(K,z))^{\frac{1}{2}}) \leq n(1,{\bf
T}(K,z)). $

In the same way one  checks that $ n(1,{\bf T} _\gamma(K,z)) \leq
n(1,2(R^0_\gamma(K,z))^{\frac{1}{2}}V
(R^0_\gamma(K,z))^{\frac{1}{2}})  . $
\end{proof}
\begin{remark}

On the left hand side of \eqref{coordinate} the operator
$V^{\frac{1}{2}} R^0_\gamma(K,z) V^{\frac{1}{2}},$ is a partial
integral operator, since the operator
 $$V^{\frac{1}{2}} f
(k_\alpha,k_3)=V^{\frac{1}{2}}_\alpha f (k_\alpha,k_3)=(I\otimes
v^{\frac{1}{2}}) f (k_\alpha,k_3) $$ is written in the coordinate
$(k_\alpha,k_3),$ that is, it is an integral operator with respect
to $k_3.$

The right hand side of \eqref{coordinate} can be written as
$V_\alpha^{\frac{1}{2}} R^0_\gamma(K,z) V_3^{\frac{1}{2}},$ where
the operator $V=V_\alpha$ is written in coordinate $(k_\alpha,k_3)$,
that is, it is integral operator with respect to $k_3.$ But the
operator $V=V_3$ is written in the coordinates $(k_3,k_\alpha)$,that
is, it is an integral operator with respect to $k_\alpha$ and hence
the operator $V^{\frac{1}{2}} R^0_\gamma(K,z) V^{\frac{1}{2}}$ on
the right hand side of \eqref{coordinate} is an integral operator in
all variables.
\end{remark}

\section{Asymptotics for the number of eigenvalues of$ H_{\gamma}(K)$}

In this section we shall prove Theorem \ref{infinite}.

\begin{theorem}\label{main} The equality
\begin{equation}\label{asimptot}
\lim\limits_{\frac{|K|^2}{M_\gamma}+|z|\to 0} \frac{n(1,{\bf
T}_{\gamma }(K,z))} {|log(\frac{|K|^2}{M_\gamma}+|z|)|}
={U}(\mu,\gamma)
\end{equation} holds.
\end{theorem}

Theorem \ref{main} will be deduced by a perturbation argument based
on Lemma 4.7, which has been proven in \cite{Sob}. For completeness,
we here reproduce the lemma.

\begin{lemma}\label{comp.pert}
 Let $A (z)=A_0 (z)+A_1 (z),$ where $A_0(z)$ (resp.$A_1(z)$) is
compact and continuous in $z<0$ (resp.$z\leq 0$).  Assume that for
some function $f(\cdot),\,\, f(z)\to 0,\,\, z\to 0-$ one has
$$
\lim_{z\rightarrow 0-}f(z)n(\lambda,A_0 (z))=l(\lambda),
$$
and $l(\lambda)$ is continuous in $\lambda>0.$ Then the same limit
exists for $A(z)$ and
$$ \lim_{z\rightarrow 0-}f(z)n(\lambda,A (z))=l(\lambda).
$$
\end{lemma}
$\Box$
\begin{remark} According
to Lemma \ref{comp.pert} any perturbation of the operator $A_0(z)$
defined in Lemma \ref{comp.pert}, which is compact and continuous up
to $z=0$ does not contribute to the asymptotics \eqref{asimptot}.
Throughout the proof of the following theorem we shall use this fact
without further comments.
\end{remark}

Let $T_{\gamma }(K,z),K\in \T^{3},z\leq \tau(K,\gamma)$ the
self-adjoint operator
 defined  in
$L_{2}(\T^{3})$ by

\begin{equation}\label{T(kz0)}
(T_{\gamma}(K,z)f)(p)=-\frac{v(\gamma)}{(2\pi)^3}\int_{\T^3}
\,\frac{\Delta_{\gamma}^{-\frac{1}{2}}
(K,p,z)\Delta_{\gamma}^{-\frac{1}{2}}(K,q,z)}
{E(K,\gamma;p,q)-z}f(q)dq.
\end{equation}

\begin{lemma}\label{eqspec}
 The equality
$$
n(1,{\bf T}_{\gamma}(K,z))=n(1,{T}_{\gamma}(K,z))
$$
holds.
\end{lemma}
\begin{proof}
Let $\Psi:L_ 2(({\T}^3)^2)\to L_ 2({\T}^3) $
 be the operator given by
$$
(\Psi f)(p)={(2\pi)^{-\frac{3}{2}}} \int_{\T^3}f(p,q)dq
$$
and let $\Psi^*$ be its adjoint.

One can easily check that the equalities
\begin{equation} \label{izom1}
\Psi f={(2\pi)^{\frac{3}{2}}} V^\frac{1}{2} f \quad\text{and}\quad
V^\frac{1}{2} W^\frac{1}{2}_\gamma
 f= \Delta^{-\frac{1}{2}}_{\gamma}(K,p,z)
V^\frac{1}{2} f,\quad f \in L_ 2(({\T}^3)^2)
\end{equation}
 hold.

These equalities  imply the equality ${\bf
T}_\gamma(K,z)=\Psi^*T_{\gamma}(K,z)\Psi.$

Since any nonzero eigenvalue of $\Psi^*T_{\gamma}(K,z)\Psi$ is an
eigenvalue of $\Psi\Psi^*T_{\gamma}(K,z)$
 as well with the same
algebraic and geometric multiplicities, and
$\Psi\Psi^{*}=I_{L_2(\T^3)},$  we have
$$
n(1,{\bf T}_\gamma(K,z))=n(1,T_{\gamma}(K,z)).
$$

\end{proof}

 \begin{lemma}\label{EK.det}
 There exists $\delta >0$ such that
 \begin{align}\label{eps exp}
E(0,\gamma;p,q)=\frac{1}{2}((1+\gamma)\,p^{2}+2\gamma\,
  (p,q)+(1+\gamma)\,q^{2} )+
 O(|p|^{4}+|q|^{4})
\end{align}
as  $p,q\rightarrow 0$ and for all $z\in(-\delta,0]$
\begin{align}\label{det exp}
\Delta_{\gamma}(0,p,z)
=\frac{v(\gamma)}{2\pi(1+\gamma)^{\frac{3}{2}}}
(n\,p^{2}-2z)^{\frac{1}{2}}+O(|p|^{2}+|z|) \quad \mbox{as} \quad
p,z\rightarrow 0,
\end{align}
where $n_\gamma=(1+2\gamma)(1+\gamma)^{-1}.$
\end{lemma}

\begin{proof}
 The asymptotics
\begin{equation}\label{eps. exp}
  \varepsilon(p)=\frac{1}{2}p^{2}+O(|p|^{4}) \quad as \quad
  p\rightarrow 0
\end{equation}
of the function $\varepsilon(p)$ yields \eqref{eps exp}.
 The
definition of $m(k,\gamma)$ and the representation \eqref{Eraz2}
gives the asymptotics
\begin{equation}\label{mkraz}
  m(k,\gamma)=\frac{\gamma}{2(1+\gamma)}k^{2} +O(|k|^{4}) \quad
  as \quad k\rightarrow 0,
\end{equation}
which yields  \eqref{det exp}.
\end{proof}

Denote by $ \chi_\delta(\cdot)$ the  characteristic function of $
U_\delta(0)=\{ p\in \T^3:\,\, |p|<\delta \}.$

Let $T(\delta,\frac{K^2}{2M_\gamma}+|z|)$ be operator on
$L_2({\T}^3)$ with the kernel
\begin{align*}
&-D_{\gamma}  \frac{ \chi_\delta (p) \chi_\delta (q) (n_\gamma p^2+
2(\frac{K^2}{2M_\gamma}+|z|))^{-1/4} (n_\gamma q^2+
2(\frac{K^2}{2M_\gamma}+|z|))^ {-1/4} }
{(1+\gamma)q^{2}+2\gamma(p,q)+(1+\gamma)p^{2}+2(\frac{K^{2}}{2M_\gamma}
+|z|)},
\end{align*}
where \[
  D_{\gamma}=\frac{(1+\gamma)^{\frac{3}{2}}}{2\pi ^{2}},\quad
 \quad n_{\gamma}=\frac{1+2\gamma}{1+\gamma} \quad
  M_\gamma =\frac{1+2\gamma}{\gamma}.
\]
\begin{lemma} \label{raznost}  The operator $ T_\gamma (K,z)-T_\gamma (\delta;
\frac{K^2}{2M}+|z|)$ belongs to the Hilbert-Schmidt class and is
continuous in $K\in \T^3$ and $z\leq 0.$
\end{lemma}
\begin{proof}

Applying the asymptotics \eqref{eps exp} and \eqref{det exp} one can
estimate the kernel of the operator $T_\gamma (K,z) -T_\gamma
(\delta; \frac{K^2}{2M_\gamma }+|z|)$ by
\begin{equation*}
 C [ (p^2+q^2)^{-1} +
|p|^{-\frac{1}{2}}(p^2+q^2)^{-1} +
(|q|^{-\frac{1}{2}}(p^2+q^2)^{-1}+1 ]
\end{equation*}
 and hence
the operator $ T_{\gamma} (K,z)-T_\gamma (\delta;
\frac{K^2}{2M}+|z|)$
 belongs to the Hilbert-Schmidt class for all $K\in U_{\delta}(0)$ and
$z \leq 0.$ In combination with the continuity of the kernel of the
operator in $K\in U_\delta (0)$ and $z<0$ this  gives   the
continuity of $T_{\gamma}
(K,z)-T_{\gamma}(\delta;\frac{K^2}{2M}+|z|)$ in $K\in U_\delta (0)$
and $z\leq 0.$
\end{proof}

Let
\begin{equation*}\label{S kernel}{\bf S}_\gamma({\bf r}):L_2((0,{\bf
r}), {\sigma_0})\to L_2((0,{\bf r}),{\sigma_0}),\,{\bf r}=1/2 | \log
(\frac{|K|^2}{2M}+|z|),\, {\sigma_0}=L_2(\S^2),\,\,
\end{equation*}
$\S^2-$ being the unit sphere in $\R^3$, be the integral operator
with the kernel
\begin{align}\label{Sobolov}
&  S_{\gamma}(t;y)=(2\pi)^{-2}\frac{u_\gamma }{\cos h y+s_\gamma t},\\
 & \,\,u_\gamma=\frac{1+\gamma}{\sqrt{1+2\gamma}},
  \,\,s_\gamma=\frac{\gamma}{1+\gamma},\\
 &y=x-x',\,x,x'\in (0,{\bf r}),\quad t=<\xi,
\eta>,\,\xi, \eta \in \S^2,\no
\end{align}

 and let
$$\hat{\bf S}_\gamma(\lambda):\,\,\sigma_ 0\rightarrow
\sigma_0,\,\,\,\lambda\in (-\infty,+\infty) $$ be the integral
operator with the following kernel

\begin{equation}\label{Stlam}
\hat {
S}_\gamma(t;\lambda)=\int\limits_{-\infty}^{+\infty}\exp{\{-i\lambda
r\}}{S}_\gamma(t;r)dr=-(2\pi)^{-1}u_{\gamma}\frac{\sinh[\lambda(arc\cos
s_{\gamma}t)]} {(1-s_{\gamma}^2t^2)^{\frac{1}{2}}\sinh
(\pi\lambda)}.
\end{equation}

 For $\mu>0,$ define
\begin{equation}\label{sobU}
 {U}(\mu,\gamma)= (4\pi)^{-1}
\int\limits_{-\infty}^{+\infty} n(\mu,\hat{\bf S}_\gamma(y))dy.
\end{equation}
\begin{lemma}\label{sobol}
 The function  $U(\mu;\gamma)$ is continuous in $\mu>0$,
 the following limit $$ \lim\limits_{{\bf r}\to \infty} \frac{1}{2}{\bf
r}^{-1}n(\mu,{\bf S_\gamma}({\bf r}))={U}(\mu;\gamma)$$ exists and
$U(\gamma)=U(1;\gamma)>0.$
\end{lemma}
\begin{remark}
 This lemma can be proven quite similarly
to the corresponding results of \cite{Sob}. In particular, the
continuity of ${U}(\mu;\gamma)$ in $\mu>0$ is a result of Lemma 3.2,
Theorem 4.5  states  the existence of the limit
$$ \lim\limits_{{\bf r}\to
\infty} \frac{1}{2}{\bf r}^{-1}n(\mu,{\bf S}_\gamma ({\bf
r})={U}(\mu,\gamma)$$ and the inequality $U(\gamma)>0$ follows from
Lemma 3.2. \end{remark}
\begin{lemma}\label{main1} The equality
\begin{equation}\label{asimptot}
\lim\limits_{\frac{|K|^2}{M_\gamma}+|z|\to 0}
\frac{n(1,T_{\gamma}(\delta,\frac{K^2}{2M_\gamma}+|z|))}
{|log(\frac{|K|^2}{M_\gamma}+|z|)|} ={U}(\gamma) \end{equation}
holds.
\end{lemma}

\begin{proof}
 The space of functions  having support in $U_\delta(0)$
 is an invariant subspace
for the operator $T_{\gamma}(\delta,\frac{K^2}{2M_\gamma}+|z|).$

Let $T_{\gamma}^{(0)}(\delta ,\frac{K^{2}}{2M_{\gamma}}+|z|)$ be the
restriction of the operator $T_{\gamma}(\delta
,\frac{K^{2}}{2M_{\gamma}}+|z|)$ on the invariant subspace
$L_{2}(U_{\delta }(0)).$

The operator \\ $T_{\gamma}^{(0)}(\delta
,\frac{K^{2}}{2M_{\gamma}}+|z|)$ is unitarily equivalent with the
operator $T_{\gamma}^{(1)}(\delta ,\frac{K^{2}}{2M_{\gamma}}+|z|)$
 acting in $L_{2}(U_{r}(0))$ by
\begin{align*}
&T_{\gamma }^{(1)}({\delta };\frac{K^{2}}{2M_\gamma}+|z|)w(p)=\\
& -D_{\gamma }\int\limits_{U_{r}(0)}\frac{(n_\gamma
p^{2}+2)^{-1/4}(n_{\gamma }q^{2}+2)^{-1/4}} {(1+\gamma
)p^{2}+2\gamma (p,q)+(1+\gamma) q^{2}+2}w(q)dq, &
\end{align*}
where $B_{r}=\{p \in T^{3}:|p|<r,\quad r=(\frac{|K|^{2}}{2M_\gamma}
+|z|)^{-\frac{1}{2}}\}.$

The equivalence is performed by the unitary dilation $$U_{r}:\,L_{2}(U_{%
\delta }(0))\rightarrow L_{2}(B_{r}),(U_{r}f)(p)=(\frac{r}{\delta }%
)^{-3/2}f(\frac{\delta }{r}p).$$

Denote by $ \chi_1(\cdot)$ the  characteristic function of $
U_1(0).$ We may replace
 $$
(n_{\gamma}p^{2}+2)^{-1/4},\,(n_{\gamma }q^{2}+2)^{-1/4}$$ and
$$(1+\gamma )p^{2}+2\gamma(p,q)+(1+\gamma )q^{2}+2$$
 by
$$(n_{\gamma}p^{2})^{-1/4}(1-\chi _{1}(p)),\quad (n_{\gamma}q^{2})^{-1/4}(1-\chi _{1}(q))$$
and $$(1+\gamma )p^{2}+2\gamma (p,q)+(1+\gamma )q^{2},$$
respectively, since the error will be a Hilbert-Schmidt operator
continuous up to\\ $\frac{|K|^{2}}{2M_\gamma}+|z|=0.$

Then we get the operator $T^{(2)}_\gamma(r)$ in $L_2 (U_r(0)
\setminus U_1(0))$ with  the kernel
\[
-D_{\gamma }\sqrt{\frac{1+\gamma  }{\sqrt{1+2\gamma }}}\frac{
|p|^{-1/2}|q|^{-1/2}}{(1+\gamma )p^{2}+2\gamma  (p,q)+(1+\gamma
)q^{2}}.
\]
By the dilation $${\bf M}:L_2(U_r(0) \setminus U_1(0))
\longrightarrow L_2((0,{\bf r})\times {\sigma_0}),\quad r=1/2|\log
(\frac{|K|^{2}}{2M_\gamma}+|z|)|$$
 where
$(M\,f)(x,w)=e^{3x/2}f(e^{ x}w),\, x\in (0,{\bf r}),\, w \in
{\S}^2,$ one sees that the operator $T_\gamma^{(2)}(r)$ is unitary
equivalent to the integral operator ${\bf S}_\gamma({\bf r}).$ The
difference of the operators ${\bf S}_\gamma({\bf r})$ and
$T_{\gamma}(\delta,\frac{K^2}{2M_\gamma}+|z|)$ is compact (up to
unitarily equivalence) and hence Lemma \ref{sobol} yields the
equality \ref{asimptot}.
 Theorem \ref{main} is thus proved.
\end{proof}

\begin{lemma}\label{mespositiv}
 For any $\gamma >0$ the inequality $\ U(\gamma)>0$ holds.
\end{lemma}

\begin{proof}
It is convenient to calculate the coefficient
$U(\gamma)=U(1;\gamma)$ by means of a decomposition of the operator
$\hat{\bf {S}_\gamma}(y)$ into the orthogonal sum over its invariant
subspaces.

Denote by $L_{l}\subset L_{2}(S^{2})$ the subspace of the  harmonics
of degree $ l=0,1,\cdots. $

It is clear that $L_{2}(S^{2})=\sum\limits_{l=0}^{\infty }\oplus
L_{l},\quad \dim L_{l}=2l+1.$

Let $P_{l}:\,L_{2}(S^{2})\rightarrow L_{l}$ be the orthogonal
projector onto $L_{l}.$

The kernel of $P_{l}$ is expressed via the Legendre polynomial
$P_{l}(\cdot ):$
\[
  {\ P}_{l}(\xi,\eta)=\frac{2l+1}{4\pi}P_{l}(<\xi,\eta >).
\]
The kernel of $\hat{\bf {S}_\gamma}(y)$ depends on the scalar
product $<\xi ,\eta
>$ only, so that the subspaces $L_{l}$ are invariant for
$\hat{\bf {S}_\gamma}(y)$ and
\begin{equation}\label{(5.9)}
  \hat{\bf {S}_\gamma}(y)=\sum\limits_{l=0}^{\infty}\oplus(\hat{\bf {S}_\gamma}^{(l)}(y)\otimes{\ P}_{l}),
 \end{equation}
where $\hat{\mathbf{S}}_{\gamma}^{(l)}(y)$ is the multiplication
operator  by the number
\begin{equation}\label{(5.10)}
  \hat{S}_{\gamma}^{(l)}(y)=2\pi\int_{-1}^{1}P_{l}(t)
  \hat{S}_{\gamma}(t;y)dt
\end{equation}
in $L_{l}$ the subspace of the harmonics of degree $l,$ and $P_l(t)$
is a Legendre polynomial. Therefore
$$
  n(\mu,\hat{\bf {S}}_{\gamma}(y))=
  \sum\limits_{l=0}^{\infty}(2l+1)n(\mu,\hat{\bf {S}}^{(l)}(y),
  \quad \mu >0
$$
 It follows from (\ref{sobU}) and  (\ref{(5.10)}) that
\begin{equation}\label{mes}
{\ U}(1, \gamma)\ge \frac{1}{4\pi}\int\limits_{-\infty}^{+\infty}
n(1,{\bf \hat{\bf {S}}}_{\gamma}^{(0)}(y))dy\ge \frac{1}{4\pi} mes
\{x:\hat{ {S}}_{\gamma}^{(0)}(x)>1 \}.
\end{equation}

By \eqref{Stlam} we first calculate \ $\hat{S}_{\gamma}^{(0)}(y):$
\begin{equation}\label{s^0}
\hat{S}_{\gamma}^{(0)}(y)=\frac{u_\gamma}{\sinh (\pi \gamma )}
\int\limits_{-1}^{1}\frac{ \sinh [y(arccos(s_\gamma
(t))]}{\sqrt{1-s^{2}_{\gamma}t^{2}}}dt.
\end{equation}

Applying the equality
\begin{equation*} \frac{u_\gamma\sinh [y(\arcsin
s_{\gamma})]}{s_{\gamma}y\cosh (\frac{\pi y}{2}
)}=\frac{u_{\gamma}}{\sinh (\pi y )}\int\limits_{-1}^{1} \frac{\sinh
[y(arccos(s_{\gamma}t))]}{\sqrt{1-s^{2}_{\gamma}t^{2}}}dt,\, y\in
\R, \end{equation*}

 Since $$\frac{\sinh [y(\arcsin
s_{\gamma})]}{s_{\gamma}y }\geq 1$$

 we have
$$\max_{y}\hat{S}_{\gamma}^{(0)}(y)\geq \frac {u_\gamma}
{\cosh(\frac{\pi y}{2})} \geq u_\gamma>1$$ This together with
\eqref{mes} completes  the proof.
\end{proof}
{\bf Proof of Theorem \ref{main}.} The proof of Theorem \ref{main}
follows from Lemmas \ref{eqspec}, \ref{raznost} and \ref{main1}.

 {\bf Proof of Theorem
\ref{infinite}.} Let the conditions  of Theorem \ref{infinite} be
fulfilled. Then the proof of Theorem \ref{infinite} follows from
Theorem \ref{main} and Lemma \ref{b-s}.

\section{Finiteness of  number of eigenvalues of $H_\gamma(K)$}
Now we are going to proof the finiteness of
$N_\gamma(K,\tau_{ess}(K))$ for $K\in \T^3$. First we shall prove
that the operator $T_\gamma(K,\tau_{ess}(K))$ belongs to the
Hilbert-Schmidt class.

The point $p=0$ is the non degenerate minimum of the functions
$\varepsilon(p)$ and $z_\gamma (p)$ (see Corollary \ref{corol}) and
hence $p=0$ is the non degenerate  minimum of $Z_\gamma (0,p)$
defined by $$ Z_\gamma (K,p):=\varepsilon(p)+ z_\gamma (K-p).$$

Using the definition of $ Z_\gamma (K,p)$ for all $K \in \T^3$ one
can conclude that minimum point $p^Z_\gamma (K)\in \T^3$ of the
function $Z_\gamma (K,p),\,K \in \T^3$ is non degenerate and  the
matrix inequality
$$
B(K)=\big (\frac{\partial^2 Z_\gamma }{\partial p^{(i)} \partial
p^{(j)}}
 (K,p^Z_\gamma (K))
 \big )_{i,j=1,2,3} >0
$$
holds. Hence   the asymptotics \be\label{Z} Z_\gamma
(K,p)=\tau_\gamma (K)+ (B(K)(p-p^Z_\gamma (K)),p-p^Z_\gamma (K))
+o(|p-p^Z_\gamma (K)|^2)\,\,\mbox{as}\,\,|p-p^Z_\gamma (K)| \to 0
\ee holds, where $\tau_\gamma (K)=Z_\gamma (K,p^Z_\gamma (K)).$ From
Lemma \ref{tasvir} we conclude that for all $K\in \T^3,p\in
U_{\delta(K)}(p^Z_\gamma (K))$ the equality
\begin{equation}\label{nondeger}
\Delta_\gamma (K,p,\tau_\gamma (K)) =(Z_\gamma (K,p)-\tau_\gamma
(K)) \hat \Delta_\gamma (K,p,_\gamma (K))
 \end{equation}
holds, where $\hat \Delta_\gamma (K,p^Z_\gamma (K),\tau_\gamma
(K))\neq 0.$ Putting \eqref{Z} into \eqref{nondeger} we get  the
following

\begin{lemma}\label{detnol} For any $K \in \T^3$  there are
  positive nonzero
constants $c$ and $C$ depending on $K$ and $U_{\delta(K)}(p^Z_\gamma
(K))$ such that  for all $p\in U_{\delta(K)}(p^Z_\gamma (K))$ the
following inequalities
 \begin{equation}\label{otsenka2}
 c|p -p^Z_\gamma (K)|^2 \leq
\Delta_\gamma (K,p,\tau_\gamma (K))\leq C|p-p^Z_\gamma (K)|^2
\end{equation}
hold.
\end{lemma}
\begin{flushright}$\Box$ \end{flushright}
\begin{remark}\label{degenerate}
Let the kernel function $v(\cdot)$ of the interaction operator $v$
and the dispersion relation $\varepsilon$ be real-analytic functions
on the three-dimensional torus $\T^3.$ In this case the minimum
(critical) values of the function $Z_\gamma(K,\cdot)$ may
degenerates only at a finite number of co-dimension $1$ manifolds
$\aleph_{n}\in \T^3,n=1,2,...,N.$
\end{remark}

\begin {lemma}\label{G-S}

 For any $K\in
\T^3$ the operator $T_{\gamma}(K, \tau_{ess}(K))$ belongs to the
Hilbert-Schmidt class.
\end{lemma}
\begin{proof}
 By Lemma \ref{tau<E} we have
\begin{equation}\label{inequal 2} \tau_{ess}(K)
 <
E_{\min}(K),\,K \in \T^3.
\end{equation}
The operator $\hbar_{\gamma}(0)$ has a zero energy resonance. By
Theorem \ref{mavjud} the operator $\hbar_{\gamma}(k),\,k\in \T^3,\,
k\neq 0$ has a unique eigenvalue $z_\gamma (k),$ $z_\gamma
(k)<m(k,\gamma).$

Since $ \tau^\gamma(K)= \min_{p \in {\T}^3}Z_\gamma(K,p)$ the
function $Z_\gamma (K,p)$ has a unique minimum and hence for all $
p\in \T^3\setminus U_\delta(p^Z_\gamma (K))$ we obtain
\begin{equation}\label{otsenka}
\Delta_\gamma(K,p,\tau^\gamma(K))\geq C>0.
\end{equation}
According to \eqref{inequal 2} for all $p_\gamma ,p_\beta\in {\T}^3$
and $K\in \T^3$ the inequality
\begin{equation}\label{inequal3}
E_{\alpha\beta}(K;p_\gamma ,p_\beta) -\tau_\gamma(K)\geq
E_{\min}(K)-\tau_{\gamma}(K)>0
\end{equation} holds.
Using \eqref{otsenka2}, \eqref{otsenka} and taking into account
\eqref{inequal3} we can make certain that for all $K\in \T^3$ and
$p_\gamma  \in U_\delta(p^Z_\gamma (K)),\, p_\beta \in
U_\delta(p^Z_\beta(K))$ the modules of the kernels $T_\gamma
(K,\tau_{ess}(K);p_\gamma ,p_\beta)$ of the integral operators
$T_\gamma (K,\tau_{ess} (K))$ can be estimated by
$$
\frac{C_0(K)}{|p_\gamma -p^Z_\gamma (K)||p_\beta-p^Z_\beta(K)|}+C_1,
$$ where $C_0(K)$ and $C_1$ some constants.
Taking into account \eqref{inequal 2} we conclude that
$$T_{\gamma}(K,\tau_{ess}(K)),\,$$ are
Hilbert-Schmidt operators.
\end{proof}
Now we shall prove the finiteness of $N_\gamma(K, \tau_{ess}(K))$
(Theorem \ref{finite}).

\begin{theorem}\label{g.b-s}
For the number $N_\gamma(K, \tau_{ess}(K))$ the relation
$$
 N_\gamma(K, \tau_{ess}(K))\leq
 \lim_{\nu \to 0}n(1-\nu,T_{\gamma}(K,\tau_{ess}(K))))
$$
holds.
\end{theorem}

 \begin{proof}  By Lemmas \ref{b-s} and \ref{eqspec} we have
$$
N_\gamma(K,z)=n(1,T_{\gamma}(K,z))\,\,\mbox{as}\,\,z<\tau_{ess}(K)
$$
and by Lemma \ref{G-S} for any $\nu\in [0,1)$ the number
$n(1-\nu,T_{\gamma}(K,\tau_{ess}(K)))),\,K \in \T^3 $ is finite.
Then according to the Weyl inequality
$$n(\lambda_1+\lambda_2,A_1+A_2)\leq
n(\lambda_1,A_1)+n(\lambda_2,A_2)$$ for all $z<\tau_{ess}(K)$ and
$\nu \in (0,1)$ we have
$$
N_\gamma(K,z)=n(1,T_{\gamma}(K,z))\leq
n(1-\nu,T_{\gamma}(K,\tau_{ess}(K))))+n(\nu,T_{\gamma}(K,z)-
T_{\gamma}(K,\tau_{ess}(K)))).
$$
Since $T_{\gamma}(K,\tau_{ess}(K))$ is continuous from the left up
to $z=\tau_{ess}(K),\,K\in \T^3$, we obtain
$$
\lim_{z\to \tau_{ess}(K)} N(K,z)= N(K,\tau_{ess}(K))\leq
n(1-\nu,T_{\gamma}(K,\tau_{ess}(K))))\,\, \mbox{for all}\,\, \nu \in
(0,1)
$$
and so $$T_{\gamma}(K,\tau_{ess}(K)))\leq \lim _{\nu\to
0}n(1-\nu,T_{\gamma}(K,\tau_{ess}(K)))).$$
\end{proof}

 Acknowledgements

This work was supported by DFG 436 USB 113/6 project. The last named
author gratefully acknowledged the hospitality of the Institute of
Applied Mathematics of the University Bonn,  the SISSA (Trieste,
Italy) and the  University Roma 1(Italy).

\end{document}